\input amstex
\documentstyle{amsppt}
\def\ds{\displaystyle}
\topmatter
\abstract{For $w\in A_p(RH_r)$ we determine the precise range of indices so 
that $w\in RH_r(A_p)$, the precise range of $q<p$ for which $w\in A_q$, and
the precise range of $\tau>1$ for which $w^\tau\in A_p$.}\endabstract
\title {  The precise range of indices for the RH$_r$- and A$_p$-
weight classes}
\endtitle
\author {C.J. Neugebauer }
\endauthor
\endtopmatter
\document
 {\bf 1. Introduction.} The weight classes which we will study in this paper
control weighted norm inequalities for the Hardy-Littlewood maximal operator,
singular integral operators, etc., and are defined as follows.  A weight $w$
, i.e., $w:\Bbb R^n\to\Bbb R^+$, is in the reverse H\"older class $RH_r, 1<r<\infty,$ provided
$\ds{\left(\frac{1}{|Q|}\int_Qw^r\right)^{1/r}\leq c\frac{1}{|Q|}\int_Qw}$, and
$w\in A_p$ iff $\ds{\frac{1}{|Q|}\int_Qw\left(\frac{1}{|Q|}\int_Qw^{1-p'}\right)^{p-1}\leq c}$ (see \cite{2}), where $Q$ is an arbitrary cube in $\Bbb R^n$. These classes are related, and we have $A_\infty
\equiv\cup_{p<\infty} A_p=\cup_{1<r} RH_r$. The purpose of this paper is the 
following.

1. For $w\in A_p$ find the precise range of $r$'s such that $w\in RH_r$, the
precise range of $q<p$ for which $w\in A_q$, and the precise range of $\tau>1$
such that $w^\tau\in A_p$.

2. For $w\in RH_r$ find the precise range of $p$'s such that $w\in A_p$, and the
precise range of $q>r$ for which $w\in RH_q$.

The higher integrability of $w\in RH_r$ is due to Gehring \cite{3}, i.e., there
exists $p>r$ such that $w\in RH_p$. Closely related to this is the property
that $w\in A_p$ implies the existence of $q<p$ such that $w\in A_q$. The
exact range of higher order integrability for $w\in RH_r$ has been found
recently by Kinnunen \cite{4} for $n=1$. The sharp upper index is given as the 
solution to $\ds{c\frac{p-r}{r}(p')^r=1}$, where $c$ is the constant in
$\ds{\frac{1}{|I|}\int_Iw^r\leq c\left(\frac{1}{|I|}\int_Iw\right)^r}$. Our
upper bound has a different form and depends upon appropriate factorizations
of $w$.  Apart from factorizations our proofs will be along the lines of
\cite{4}. 

In sections 2 and 3 our analysis will be in $\Bbb R$, and $\Bbb R^n$-versions
of these results will be taken up in section 4. 
\vskip.3truein

{\bf 2. The classes $\bold{A_1}$ and $\bold{RH_\infty}$.} 
Throughout $I$ will denote an
arbitrary interval in $\Bbb R$. For $w\in A_1$, we let $A_1(w)$ be the $\inf$
of all constants $c$ for which $\ds{\frac{1}{|I|}\int_Iw\leq c\inf_Iw}$, and
for $w\in RH_\infty$, we denote by $RH_\infty(w)$ the $\inf$ of all constants
$c$ such that $\ds{\sup_Iw\leq c\frac{1}{|I|}\int_Iw}$ for all intervals $I$.
The classes $A_1$ and $RH_\infty$ generate $A_\infty$ in the sense that $A_\infty=A_1\cdot
RH_\infty$, and $RH_\infty$ plays the same role relative to $RH_r$ as $A_1$
does to $A_p$ \cite{1}.

\proclaim{Theorem 1} Let $w\in A_1$ and let $c=A_1(w)$. Then $w\in RH_r$ for
all $\ds{1<r<\frac{c}{c-1}}$. The range of $r$'s is best possible on $\Bbb R_+$.
\endproclaim

\demo{Proof} This is one of the results in Kinnunen's thesis \cite{4}.  The proof
consists by showing first that for $\lambda\geq\inf_Iw\equiv \lambda'_I$
$$
\int_{\{x\in I:w(x)>\lambda\}}w\,dx\leq c\lambda |\{x\in I:w(x)>\lambda\}|\tag1
$$
where $c=A_1(w)$. This is a reverse Chebyshev inequality and is in fact
equivalent with $w\in A_1$ \cite{4}. Next multiply this inequality by
$\lambda^r, r>-1$ and integrate with respect to $\lambda$ from $\lambda'_I$ to
$\infty$. We shall be more precise in section 4 when we treat the $\Bbb R^n$-
case using the same type of distributional inequality. The weight $x^{1/c-1},
c>1$ shows that the range of $r$ is best possible on $\Bbb R_+$.
\enddemo

{\smc Remark:} A. Torre pointed out to me that on $\Bbb R$ the weight
$w(x)=|x|^{1/c-1}, c>1$ does not provide an example for the precise range since
$A_1(w)\leq 2c$.  This is easily seen by considering intervals $I=[a,b]$ with
$a<0<b$ and $|a|\leq b$. However, if we define for $w\in A_1$, $c_w=
\inf c$ for which $\ds{\frac{1}{|I|}\int_Iw\leq c\inf_Iw}$ for all intervals
$I\subset \Bbb R\setminus \{0\}$, then $w\in RH_p, 1<p<c_w/(c_w-1)$.
To prove this, observe first that $w\in RH_p(\Bbb R_+)\cap RH_p(\Bbb R_-)$
for $1<p<c_w/(c_w-1)$ by Theorem 1.  If $I=[a,b]$ with $a<0<b, |a|\leq b$, and 
if $I_0=[-b,b]$, then
$$\gather
\frac{1}{|I|}\int_Iw^p\leq\frac{2}{|I_0|}\int_{I_0}w^p\leq \frac{1}{b}\int_{-b}
^0w^p+\frac{1}{b}\int_0^bw^p\\\leq c\left((\frac{1}{b}\int_{-b}^0w)^p+
(\frac{1}{b}\int_0^bw)^p\right)\leq c_p(\frac{1}{|I|}\int_Iw)^p,
\endgather
$$
for $1<p<c_w/(c_w-1)$. This range is best possible, as the weight $w(x)=
|x|^{1/c-1}, c>1,$ shows.  The same type of remark applies to the remaining
Theorems of this note.

\proclaim{Theorem 2} Let $w\in RH_\infty$ and let $c=RH_\infty(w)$. Then
$w\in A_p$ for all $p>c$, and the range of $p$'s is best possible on $\Bbb R_+$.
\endproclaim

\demo{Proof} We first show that for $0<\lambda\leq \lambda_I\equiv\sup_Iw$
$$
\lambda|\{x\in I:w(x)<\lambda\}|\leq c\int_{\{x\in I:w(x)<\lambda\}}w\tag2
$$
where $c=RH_\infty(w)$. If $G$ is an arbitrary open set with $I\cap\{w<\lambda\}\subset G\subset I$, then we have to show that $\ds{L\leq c\int_Gw}$, where
$L$ is the left side of the above inequality.

There are two cases.  If $|G|=|I|$, then $L\leq\lambda|G|=\lambda |I|\leq c
\int_Iw=\int_Gw$. If $|G|<|I|$, we cover $G$ by a disjoint union of intevals
$\{I_j\}$ maximal with respect to $|I_j\setminus G|=0$.  Then $|G|=\sum|I_j|$,
and $|\sigma I_j\setminus G|>0$ for every $\sigma>1$, where $\sigma I$ is the
interval concentric with $I$ and of length $\sigma|I|$. It follows that
$\ds{\sup_{\sigma I_j}w\geq \lambda}$ and thus for every $\sigma>1$
$$
\lambda|\sigma I_j|\leq c\int_{\sigma I_j}w.
$$
We now let $\sigma\searrow 1$ and obtain $\ds{\lambda|I_j|\leq c\int_{I_j}w}$.
Consequently
$$
L\leq\sum\lambda|I_j|\leq c\sum\int_{I_j}w=c\int_Gw.
$$

Below we need the finiteness of $\int_I1/w$ and for that purpose we truncate
$w$ from below, i.e., for $0<\beta<\lambda\leq \lambda_I$ let
$$
w_\beta(x)=\cases w(x),&w(x)>\beta\\
\beta,&w(x)\leq\beta.\endcases
$$
Since $RH_\infty(w_\beta)\leq c$ we have the same inequality (2) with $w$ replaced by $w_\beta$.  We multiply this inequality by $\lambda^{-r}, r>2$ and integrate to get
$$
\int^{\lambda_I}_\beta\frac{1}{\lambda^{r-1}}\int\chi_{\{w_\beta<\lambda\}\cap I}(x)\,dx\,d\lambda\leq c\int^\infty_\beta\frac{1}{\lambda^r}\int_{\{w_\beta<
\lambda\}\cap I}w\,dx\,d\lambda.
$$
In the left side of the above inequality we interchange the order of integration and obtain
$$\gather
\int_I\int^{\lambda_I}_{w_\beta(x)}\frac{1}{\lambda^{r-1}}\,d\lambda\,dx=
\frac{1}{r-2}\int_I\left(\frac{1}{w_\beta^{r-2}}-\frac{1}{\lambda_I^{r-2}}
\right)\\
=\frac{1}{r-2}\int_I\frac{dx}{w_\beta^{r-2}}-\frac{1}{r-2}\frac{|I|}{\lambda
_I^{r-2}}.\endgather
$$
Similarly, the right side equals
$$\int_I\int^\infty_{w_\beta(x)}w(x)\frac{1}{\lambda^r}\,d\lambda\,dx=
\frac{1}{r-1}\int_I\frac{dx}{w_\beta^{r-2}}.
$$
Hence
$$
\int_I\frac{dx}{w_\beta^{r-2}}\leq c\frac{r-2}{r-1}\int_I\frac{dx}{w_\beta^{r-2}}+\frac{|I|}{\lambda_I^{r-2}}.
$$
Next choose $r>2$ so that $c(r-2)<r-1$.  Since the integrals involved are $<\infty$ we see that
$$
C\int_I\frac{dx}{w_\beta^\alpha}\leq\frac{|I|}{\lambda_I^\alpha},
$$
where $\alpha=r-2$.  We let $\beta\searrow 0$ and get that $w^{-\alpha}\in A_1$.
This, of course, immediately implies that $w\in A_{1+1/\alpha}$. Also observe
that $c\alpha<\alpha+1$ is equivalent with $p=1+1/\alpha>c$.

The weight $w(x)=x^{c-1},\quad c>1$ has $RH_\infty(w)=c$ on $\Bbb R_+$ and is in $A_p$
for $p>c$, but not in $A_c$.
\enddemo

{\smc Remark:}(i) It is easy to see that (2) is actually equivalent with $w\in
RH_\infty$. (ii) For an example of the best range on $\Bbb R$ proceed as in
the remark after Theorem 1.
\vskip.3truein
{\bf 3. The classes $\bold{A_p}$ and $\bold{RH_r}$}. In this section we use 
Theorems 1,2 for the extension to the classes $A_p$ and $RH_r$, and, as we
shall see, the range of the indices will be governed by factorizations of
the weight.  In \cite{1} we have shown that $w\in RH_r$ iff $w=w_0w_1$ where
$w_0\in RH_\infty$ and $w_1\in RH_r\cap A_1$. But then $w_1^r\in A_1$, and thus
$w\in RH_r$ iff $w=uv^{1/r}$ with $u\in RH_\infty$ and $v\in A_1$.  We shall
also use the factorization \cite{2} for $A_p$, i.e., $w\in A_p$ iff $w=uv^{1-p}$
with $u,v\in A_1$.

\proclaim{Theorem 3} Let $w=uv^{1/r}$ be in $RH_r$ with $u\in RH_\infty$ and
$v\in A_1$, and let $c_1=RH_\infty(u)$.  Then $w\in A_p$ for all $p>c_1$. This
range of $p$'s is best possible on $\Bbb R_+$.\endproclaim

\demo{Proof} Let $p>c_1$.  Then
$$\gather
\frac{1}{|I|}\int_Iuv^{1/r}\left(\frac{1}{|I|}\int_I(uv^{1/r})^{1-p'}\right)^{p-1}\leq\\
\sup_Iu\frac{1}{|I|}\int_Iv^{1/r}\left(\frac{1}{|I|}\int_Iu^{1-p'}\right)^{p-1}
\sup_I\frac{1}{v^{1/r}}\leq\\c\frac{1}{|I|}\int_Iu\inf_Iv^{1/r}\left(\frac{1}{|I|}\int_Iu^{1-p'}\right)^{p-1}\sup_I\frac{1}{v^{1/r}}\leq C,
\endgather
$$
since $u\in A_p, p>c_1$ by Theorem 2.

The weight $w(x)=x^{c_1-1}$ which is in $RH_\infty\subset RH_r$
shows that the range of $p$'s is best possible on $\Bbb R_+$.
\enddemo

The next result will give us the precise range of higher integrability of
$w\in RH_r$.

\proclaim{Theorem 4}  Let $w=uv^{1/r}$ be in $RH_r$  with $u\in RH_\infty$ and
$v\in A_1$.  If $c_2=A_1(v)$, then $w\in RH_p$ for all $r\leq p<c_2r/(c_2-1)$.
The range of $p$'s is best possible on $\Bbb R_+$.\endproclaim

\demo{Proof} Let $p$ satisfy the above inequality, and then choose $q>1$ such 
that
$$p<\frac{c_2r}{q(c_2-1)}.$$
Since $\ds{1\leq \frac{pq}{r}<\frac{c_2}{c_2-1}}$, by Theorem 1, $v\in RH_{qp/r}$.  This and H\"older's inequality
gives us
$$\gather
\frac{1}{|I|}\int_Iw^p=\frac{1}{|I|}\int_Iu^pv^{p/r}\leq\left(\frac{1}{|I|}
\int_Iu^{q'p}\right)^{1/q'}\left(\frac{1}{|I|}\int_Iv^{qp/r}\right)^{1/q}\\
\leq \sup_Iu^p\cdot c\left(\frac{1}{|I|}\int_Iv\right)^{p/r}\leq c'\sup_Iu^p\cdot \inf_Iv^{p/r}\leq c''\left(\frac{1}{|I|}\int_Iuv^{1/r}\right)^p,\endgather
$$
since $u\in RH_\infty$.

Let $0<\alpha<1$ and consider the weight $w(x)=x^{-\alpha}$ on $\Bbb R_+$..  Then $w\in
RH_r$ for $1<r<1/\alpha$.  We fix such an $r$ and write $w=v^{1/r},\quad
v=x^{-\alpha r}$. This is the factorization $w=uv^{1/r}$ with $u\equiv 1$.
Since $A_1(v)=c_2=1/(1-\alpha r)$ we see that $c_2r/(c_2-1)=1/\alpha$ which
is the precise upper bound of higher integrability for this weight.
\enddemo

For the next Theorem we need the fact \cite{1} that $v\in A_1$ implies that
$(1/v)^\gamma\in RH_\infty$ for every $\gamma>0$.

\proclaim{Theorem 5}  Let $w=uv^{1-p}$ be in $A_p$ with $u,v\in A_1$, and let
$c=A_1(u)$.  Then $w\in RH_r$ for all $1<r<c/(c-1)$.  The range
of $r$'s is best possible on $\Bbb R_+$.\endproclaim

\demo{Proof} We use Theorem 1 and observe that
$$\gather
\frac{1}{|I|}\int_Iw^r=\frac{1}{|I|}\int_Iu^r\frac{1}{v^{r(p-1)}}\leq\\
\sup_I\frac{1}{v^{r(p-1)}}\frac{1}{|I|}\int_Iu^r\leq C\sup_I\frac{1}{v^{r(p-1)}}\left(\frac{1}{|I|}\int_Iu\right)^r\leq\\
C\,A_1(u)^r\sup_I\frac{1}{v^{r(p-1)}}\cdot(\inf_Iu)^r\leq\\
C\,A_1(u)^r\cdot c^r\left(\frac{1}{|I|}\int_I\frac{1}{v^{p-1}}\cdot\inf_Iu\right)^r\leq\\
C\,A_1(u)^rc^r\left(\frac{1}{|I|}\int_Iw\right)^r.\endgather
$$

The example in Theorem 1 shows that the range of $r$'s is best possible on $\Bbb R_+$..
\enddemo

\proclaim{Corollary.} Let $w=uv^{1-p}$ be in $A_p$ with $u,v\in A_1$.  If
$c=\max\{A_1(u),A_1(v)\}$, then $w^\tau\in A_p$ for $1\leq\tau<c/(c-1)$ and
the range of $\tau$'s is best possible on $\Bbb R_+$.\endproclaim

\demo{Proof} From Theorem 5 we have that $w$ and $w^{1-p'}$ are in $RH_\tau$
for $1\leq\tau<c/(c-1)$ and hence $w^\tau\in A_p$.  The weight $w(x)=
x^{-\alpha}, 0<\alpha<1$ is in $A_1\subset A_p$ on $\Bbb R_+$, and thus we can take
$u=x^{-\alpha}, v\equiv 1$.  Then $c=1/(1-\alpha)$ and thus $c/(c-1)=1/\alpha$, which is best possible.\enddemo

We will now use Theorem 5 to get the exact range on $q<p$ such that $w\in A_p$ implies
$w\in A_q$.

\proclaim{Theorem 6}  Let $w=uv^{1-p}$ be in $A_p$ with $u,v\in A_1$ and let
$c_*=A_1(v)$.  Then $w\in A_q$ for all $q$ satisfying
$$
\frac{(p-1)(c_*-1)}{c_*}+1<q\leq p.
$$
The range of $q$'s is best possible on $\Bbb R_+$.\endproclaim

\demo{Proof}  Since $w^{1-p'}=vu^{1-p'}$ is in $A_{p'}$ we get from Theorem 5 that $w^{1-p'}
\in RH_r$ for $1<r<c_*/(c_*-1)$. Hence since $w\in A_p$
$$
\frac{1}{|I|}\int_Iw\left(\frac{1}{|I|}\int_Iw^{r(1-p')}\right)^{(p-1)/r}\leq
c\frac{1}{|I|}\int_Iw\left(\frac{1}{|I|}\int_Iw^{1-p'}\right)^{p-1}\leq C<\infty.
$$
Thus $w\in A_q,\quad q=1+(p-1)/r$, i.e.,$(p-1)(c_*-1)/c_*+1<q\leq p.$

We will now show that this range is best possible.  Let on $\Bbb R_+$, $w(x)=x$ and fix
$p_0>2$.  Then $w\in A_{p_0}$ and $w=v^{1-p_0}$ with $v=x^{1-p_0'}\in A_1$.
Since $A_1(v)=1/(2-p_0')=c_*$ the lower bound of the range of $q$'s given above
is $2$. It is known \cite{2} that $w\in A_q$ exactly for $q>2$.\enddemo
\vskip.3truein

{\bf 4. Extensions to $\Bbb R^n$.} The results in $\Bbb R^n$, $n>1$, are
somewhat different from the $n=1$ versions, and the reason is that we do not 
know whether the covering of an open set used in Theorem 2 , i.e., $G\subset
\cup I_j$, where the $I_j$'s are disjoint, $|I_j\setminus G|=0$, and for
$\sigma>1,\quad |\sigma I_j\setminus G|>0$, has a corresponding analogue in
$\Bbb R^n$. To avoid this difficulty, we consider for $k\geq 3$ the classes
$A_{1,k}$ and $RH_{\infty,k}$ defined as follows.  We say that $w\in A_{1,k}$
iff $\ds{\frac{1}{|Q|}\int_{Q}w\leq c'_k\inf_{kQ}w}$, and we say that
$w\in RH_{\infty,k}$ iff $\ds{\sup_{kQ}w\leq c''_k\frac{1}{|Q|}\int_Qw}$, where
$Q$ is the generic notation of a cube in $\Bbb R^n$ and $kQ$ is the cube concentric with $Q$ having side-length $k\times$ the side-length of $Q$. 

 It is easily seen that $A_{1,k}=A_1$ and $RH_{\infty,k}=RH_\infty$.  The first
 equality is obvious, and in the second use the fact that $w\in RH_\infty$ is
doubling. It is also clear that $A_1(w)\leq c'_k$ and $RH_\infty(w)\leq c''_k$.

\proclaim{Theorem 7} Let $w\in A_{1,k}$ with constant $c'_k$. Then $w\in RH_r$
for $1<r<c'_k/(c'_k-1)$.\endproclaim

\demo{Proof} This is the n-dimensional version of Theorem 1, and we follow the
proof in \cite{4}.  We first establish
$$
\int_{\{x\in Q:w(x)>\lambda\}}w\leq c'_k\lambda|\{x\in Q:w(x)>\lambda\}|\tag3
$$
for every $\lambda\geq \inf_Qw\equiv \lambda'_Q$. To prove this, let $G$ be
an arbitrary open set with $\{w>\lambda\}\cap Q\subset G\subset Q$.  We need
to show that
$$
\int_Gw\leq c'_k\lambda |G|.
$$
This follows if $|Q|=|G|$.  If $|G|<|Q|$, we let $\Delta$ be the collection
of the dyadic cubes in $Q$ generated by $Q$.  For each $x\in G$, let 
$\tilde Q_1(x)\in\Delta$ be the cube of least side-lenght containing $x$
such that $|\tilde Q_1\setminus G|>0$. The next generation dyadic subcube of 
$\tilde Q_1$ containing $x$, say $Q_1(x)$ has the property that $|Q_1(x)
\setminus G|=0$. Since we are dealing with dyadic cubes, we can in this way
write $G\subset\cup Q_j$, where the $Q_j$'s are non-overlapping, $|Q_j\setminus
 G|=0$, and $|kQ_j\setminus G|>0$ for each $j$.  Thus $\inf_{kQ_j}w\leq \lambda$.  From this we see that
$$
\int_Gw\leq\sum\int_{Q_j}w\leq c'_k\lambda\sum |Q_j|=c'_k\lambda |G|.
$$
The rest of the argument is exactly the same as in \cite{4}.
\enddemo

{\smc Remark:} The best range of Theorem 7 is not as exact as the range of
Theorem 1, but it is so within $\epsilon>0$. We work on $\Bbb R_+$, and there
$w\in A_{1,k}$ means: For every $I\subset\Bbb R_+$,
$$
\frac{1}{|I|}\int_Iw\leq c'_k\inf_{kI\cap\Bbb R_+}w, k>1.
$$
Then by Theorem 7, $w\in RH_r(\Bbb R_+)$, $1<r<c'_k/(c'_k-1)$.  Let
now $w(x)=x^{-\gamma}, 0<\gamma<1$, and let $c_{k,\gamma}=\inf c'_k$ for
which the above displayed inequality holds.  We claim that
$$
c_{k,\gamma}=\frac{1}{1-\gamma}\left(\frac{k+1}{2}\right)^\gamma.
$$
If $I=[\alpha,\beta], \alpha\geq 0$, then the right endpoint of $kI$ is
$b=\alpha+(\beta-\alpha)(k+1)/2$.  If we let $J=[0,\beta]$, then
$$
\frac{1}{|I|}\int_Ix^{-\gamma}\leq\frac{1}{|J|}\int_Jx^{-\gamma}=\frac{1}
{1-\gamma}\left(\frac{k+1}{2}\right)^\gamma b_1^{-\gamma},
$$
where $b_1=\beta(k+1)/2=\inf_{kJ\cap\Bbb R_+}w$.  Since $k>1$, $b_1>b$ , and thus $b_1^{-\gamma}<
b^{-\gamma}$, and our claim follows.

Next we note that
$$
\frac{c_{k,\gamma}}{c_{k,\gamma}-1}=\frac{(k+1)^\gamma}{(k+1)^\gamma-
2^\gamma(1-\gamma)}<\frac{1}{\gamma}.
$$
Let now $\epsilon>0$.  Since, as $\gamma\nearrow 1$, $1/(\gamma+\epsilon)\to
1/(1+\epsilon)$, and since $c_{k,\gamma}/(c_{k,\gamma}-1)>1$, there is
$0<\gamma=\gamma_\epsilon<1$ such that
$$
\frac{1}{\gamma+\epsilon}<\frac{c_{k,\gamma}}{c_{k,\gamma}-1}<\frac{1}{\gamma}.
$$
Finally, note that $x^{-\gamma}\in RH_r(\Bbb R_+)$ for $1<r<1/\gamma$.

\proclaim{Theorem 8}  Let $w\in RH_{\infty,k}$ for some $k\geq 3$ with
constant $c''_k$.  Then $w\in A_p$ for all $p>c''_k$.  
\endproclaim

\demo{Proof} The proof is the same as the proof of Theorem 2 using the
dyadic covering of Theorem 7.
\enddemo

{\smc Remark:} Again the best range of $p'$s in Theorem 8 is not as precise as
the range in Theorem 2, but it is so within $\epsilon>0$.  As before, on
$\Bbb R_+$, $w\in RH_{\infty,k}$ means
$$
\sup_{kI\cap\Bbb R_+}w\leq c''_k\frac{1}{|I|}\int_Iw, k>1.
$$
Let $w(x)=x^r, r>0$, and let $c_{k,r}=\inf c''_k$ appearing above.  We claim
that
$$
c_{k,r}=(r+1)\left(\frac{k+1}{2}\right)^r.
$$
To see this, let $I=[\alpha,\beta], \alpha>0$.  Then the right endpoint of
$kI$ is $b=\alpha+(\beta-\alpha)(k+1)/2$.  Since $x^r$ increases,
$$
\frac{1}{|I|}\int_Iw\geq\frac{1}{\beta}\int_0^\beta w=\frac{1}{r+1}\left(
\frac{2}{k+1}\right)^rb_1^r,
$$
where $b_1=\beta(k+1)/2>b=\sup w$ on $kI$.  Let $\epsilon>0$ be given.
	Then there is $0<r=r_\epsilon$ such that $r+1+\epsilon>c_{k,r}>
r+1$.  The right inequality is obvious for any $r>0$, and the left
inequality follows from the fact that as $r\to 0$, $c_{k,r}\to 1$ and
$r+1+\epsilon\to 1+\epsilon$.  Finally note that $x^r\in A_p(\Bbb R_+)$
precisely when $p>r+1$.

{\smc Remark:} The remaining Theorems in $\Bbb R^n$ corresponding to the
Theorems 3,4,5,6 are the same with the proper change of constants.

\medskip
\noindent
{\smc Department of Mathematics, Purdue University, Lafayette, IN 47907}
\newline
{\it E-mail address:} neug\@math.purdue.edu
\enddocument